\documentstyle[12pt]{article}
\begin{document}

\newcommand{\mysection }[1]{\section{#1}\setcounter{equation}{0}}

\title{\bf Manifolds with Pointwise Ricci Pinched Curvature} \vspace{15mm}
\author{\bf}
\date{}
\maketitle \centerline{\large \bf Hui-Ling Gu } \vspace{8mm}
\centerline{{\large Department of Mathematics}} \vspace{3mm}
\centerline{{\large Sun Yat-Sen University }} \vspace{3mm}
\centerline{{\large Guangzhou, P.R.China}} \vspace{3mm}
\vspace{15mm}

\noindent {\bf Abstract } In this paper, we proved a compactness
result about Riemannian manifolds with an arbitrary pointwisely
pinched Ricci curvature tensor.

\begin{center}

{\bf \large \bf{1. Introduction}}
\end{center}

Let $M^n$ be an $n$-dimensional complete Riemannian manifold with
$n\geq 3$. One of the basic problems is under which condition on
its curvature the Riemannian manifold is compact. The classical
Bonnet-Myers' theorem states that a complete Riemannian manifold
with positive lower bound for its Ricci curvature is compact.

In \cite{Ha94}, Hamilton proved that:

 {\it Any convex hypersurface with dimension $\geq 3$ in Euclidean space
with second fundamental form $h_{ij}\geq
\delta\cdot\frac{tr(h)}{n}$ must be compact.}

In \cite{CZ}, Chen-Zhu proved an intrinsic analogue of the
Hamilton's result by using the Ricci flow which was introduced by
Hamilton in 1982. They proved that:

 {\it If $M^n$ is a complete $n$-dimensional $(n\geq 4)$ Riemannian
manifold with positive and bounded scalar curvature and satisfies
the following pointwisely pinching condition
$$|W|^2+|V|^2\leq \delta_n(1-\varepsilon)^2|U|^2,$$ for
$\varepsilon>0,\delta_4=\frac{1}{5},\delta_5=\frac{1}{10}$ and
$\delta_n=\frac{2}{(n-2)(n+1)},(n\geq 6)$, where $W,V,U$ denote
the Weyl conformal curvature tensor, traceless Ricci part and the
scalar curvature part of the curvature operator respectively. Then
$M^n$ is compact.}

For the 3-dimensional case, they weaken the curvature operator
pinching condition to an arbitrary Ricci curvature pinching
condition:

 {\it Let $M$ be a complete 3-dimensional Riemannian manifold with bounded
 and nonnegative sectional curvature. If $M$ satisfies the positive Ricci pinching condition:
 $$R_{ij}\geq \varepsilon\cdot scal \cdot g_{ij}>0$$ for some $\varepsilon>0$. Then $M$ must be compact.}

Recently, by the Ricci flow and the new invariant cone
construction introduced by B$\ddot{o}$hm-Wilking \cite{BW}, Ni-Wu
\cite{NW} proved the following compactness result in terms of
curvature operator:

{\it If $M^n$ is a complete $n$-dimensional $(n\geq 3)$ Riemannian
manifold with bounded curvature and satisfies $$Rm\geq \delta
U>0$$ for $\delta>0$, where $Rm,U$ denote the curvature operator
and its scalar curvature part. Then $M^n$ must be compact.}

Naturally, from the above results, one expects that: any complete
Riemannian manifold with dimension $\geq 3$ and has positive Ricci
pinched curvature must be compact. This is already true in
3-dimensional case by the result in Chen-Zhu \cite{CZ}. In this
paper, by using the Yamabe flow, we give an affirmative answer in
the class of locally conformally flat manifolds. Our main result
is the following:

\vskip 0.3cm \noindent {\bf Theorem 1.1} \emph{ Let $n\geq 3$.
 Suppose $M^n$ is a smooth complete locally conformally flat $n$-dimensional manifold
 with bounded and positive scalar
 curvature. Suppose $M^n$ has nonnegative sectional curvature and satisfies the following Ricci
 curvature pinching condition $$R_{ij}\geq \varepsilon\cdot scal\cdot g_{ij} \eqno(1.1)$$
 for some $\varepsilon >0$. Then $M^n$ is compact.}

We briefly describe the proof of the theorem. Our proof of Theorem
1.1 depends on the Yamabe flow and the limit solution of Yamabe
flow. Suppose there exists such a noncompact Riemannian manifold
satisfying the Ricci pinching condition (1.1), we evolve it by the
Yamabe flow. By the short-time existence result \cite{CZ02} and
the Ricci pinching condition, we can obtain a long-time existence
result. In section 2, we will study the asymptotic behaviors of
the solution to the Yamabe flow. Finally in section 3, we will
complete the proof of the main theorem by using the results
obtained in section 2.

\begin{center}{ \bf \large \bf {2. The Asymptotic Behaviors of the Yamabe Flow}}\end{center}

In the geometric flows, in order to know the initial manifold
well, we usually need to study the asymptotic behaviors of the
solution of the flow. In this section, we study the asymptotic
behaviors of the Yamabe flow. First we recall the Li-Yau-Hamilton
inequality of Chow \cite{Chow2} on locally conformally flat
manifolds.

\noindent {\bf Theorem 2.1} (Chow \cite{Chow2})  \emph{ Suppose
$(M^n,g_{ij})$ is a smooth $n$-dimensional $(n\geq 3)$ complete
locally conformally flat manifold with bounded and nonnegative
Ricci curvature. Let $R(x,t)$ be the scalar curvature of the
solution of the Yamabe flow with $g_{ij}$ as initial metric. Then
we have
$$\frac{\partial R}{\partial t}+\langle \nabla
R,X\rangle+\frac{1}{2(n-1)}R_{ij}X^iX^j+\frac{R}{t}\geq 0$$ for
any vector $X$ on $M$.}

In his paper \cite{Chow2}, Chow proved the above theorem for
compact locally conformally flat manifolds with positive Ricci
curvature. However, by a perturbation argument as in \cite{Ha93},
it is clear that the Li-Yau-Hamilton inequality actually holds for
complete locally conformally flat manifolds with nonnegative Ricci
curvature.

\noindent {\bf Lemma 2.2}  \emph{ Let $g_{ij}(t)$ be a locally
conformally flat complete solution to the Yamabe flow for $t>0$
which has bounded and positive Ricci curvature. If the Harnack
quantity $$Z=\frac{\partial R}{\partial t}+\langle \nabla
R,X\rangle+\frac{1}{2(n-1)}R_{ij}X^iX^j+\frac{R}{t}$$ is positive
for all $X\in T_{x_0}M^n$ at some point $x=x_0$ and $t=t_0>0$,
then it is positive for all $X\in T_{x}M^n$ at every point $x\in
M^n$ for any $t>t_0$. }
 \vskip 0.1cm \noindent{\bf Proof.} By the calculation in \cite{Chow2}, we
know
$$(\frac{\partial}{\partial t}-(n-1)\triangle )Z\geq
 (R-\frac{2}{t})Z\geq -\frac{2}{t}Z. \eqno (2.1)$$
Since $Z$ is positive for all $X\in T_{x_0}M^n$ at $t=t_0>0$, we
can find a nonnegative function $F$ on $M^n$ with support in a
neighborhood of $x_0$ so that $F(x_0)>0$ and $Z\geq
\frac{F}{t_0^2}$ for all $X$ everywhere at $t=t_0$. Let $F$ evolve
by the heat equation $$\frac{\partial F}{\partial
t}=(n-1)\triangle F. \eqno (2.2)$$ It then follows the usual
strong maximum principle that $F>0$ everywhere for any $t>t_0$. We
only need to prove that $$Z\geq \frac{F}{t^2},\qquad \mbox{for
all} \quad t\geq t_0.$$ By (2.1) and (2.2) we know
$$(\frac{\partial}{\partial t}-(n-1)\triangle )(Z-\frac{F}{t^2})\geq
 -\frac{2}{t}(Z-\frac{F}{t^2}), $$
for $t\geq t_0$. By the maximum principle we get $Z\geq
\frac{F}{t^2}$.

 This completes the proof of the Lemma 2.2.$$\eqno
\#$$

Before we give the main result of this section, we first recall
some definitions for the classification of the asymptotic
behaviors of the solution of the Yamabe flow as $t\rightarrow
+\infty$.

\noindent {\bf Definition 2.3}  (i) \emph{A complete solution to
the Yamabe flow is called a Type I limit solution if the solution
has nonnegative Ricci curvature and exists for $-\infty <t<\Omega$
for some constant $\Omega$ with $0<\Omega <+\infty$ and $R\leq
\frac{\Omega}{\Omega-t}$ everywhere with equality somewhere at
$t=0$}.

(ii) \emph{A complete solution to the Yamabe flow is called a Type
II limit solution if the solution has nonnegative Ricci curvature
and exists for $-\infty <t<+\infty$ and $R\leq 1$ everywhere with
equality somewhere at $t=0$}.

(iii) \emph{A complete solution to the Yamabe flow is called a
Type III limit solution if the solution has nonnegative Ricci
curvature and exists for $-A <t<+\infty$ for some constant $A$
with $0<A<+\infty$ and $R\leq \frac{A}{A+t}$ everywhere with
equality somewhere at $t=0$}.

\noindent {\bf Definition 2.4} (i) \emph{ We call a solution to
the Yamabe flow a steady soliton, if it satisfies
$$Rg_{ij}=g_{jk}\nabla_iX^k,$$ where $X^i$ is a vector field on the manifold.}

(ii) \emph{ We call a solution to the Yamabe flow a shrinking
soliton, if it satisfies
$$(R-\lambda)g_{ij}=g_{jk}\nabla_iX^k,$$ where $X^i$ is a vector field on the manifold and
$\lambda$ is a positive constant.}

(iii) \emph{ We call a solution to the Yamabe flow an expanding
soliton, if it satisfies
$$(R+\lambda)g_{ij}=g_{jk}\nabla_iX^k,$$ where $X^i$ is a vector field on the manifold and
$\lambda$ is a positive constant.}

\emph{Moreover, if the vector field $X$ is the gradient of some
function $f$, then we will call the corresponding soliton a
 steady, shrinking, expanding gradient soliton respectively.}

We now follow Hamilton \cite{Ha93E} and Chen-Zhu \cite{CZ} (or
also Cao \cite{Cao}) to give a classification for Type II and Type
III limit solutions.

\noindent {\bf Theorem 2.5} \emph{ Let $M^n$ be a smooth
$n$-dimensional locally conformally flat and simply connected
Riemannian manifold. Then:}

(i)\emph{ any Type II limit solution with positive Ricci curvature
to the Yamabe flow on $M^n$ is necessarily a homothetically steady
gradient soliton;}

(ii)\emph{ any Type III limit solution with positive Ricci
curvature to the Yamabe flow on $M^n$ is necessarily a
homothetically expanding gradient soliton.}

\vskip 0.1cm \noindent{\bf Proof.} The following arguments are
adapted from Hamilton \cite{Ha93E} and Chen-Zhu \cite{CZ} (or also
Cao \cite{Cao}), where the classification for the limit solutions
of the Ricci flow were given. We only give the complete proof of
(ii), since the proof of (i) is similar and easier. At the end of
the proof we point the difference between (i) and (ii), and then
it is easy to see that the rest of the arguments are the same.

By the definition of the Type III limit solution, after a shift of
the time variable, we may assume the Type III limit solution
$g_{ij}(t)$ is defined for $0<t<+\infty$ with uniformly bounded
curvature and positive Ricci curvature where $tR$ assumes its
maximum in space-time.

Suppose $tR$ assumes its maximum at a point $(x_0,t_0)$ in
space-time, then $t_0>0$ and the Harnack quantity
$$Z=\frac{\partial R}{\partial t}+\langle \nabla
R,X\rangle+\frac{1}{2(n-1)}R_{ij}X^iX^j+\frac{R}{t},\eqno (2.3)$$
vanishes in the direction $X=0$ at $(x_0,t_0)$. By Lemma 2.2 we
know that at any earlier time $t<t_0$ and at every point $x\in
M^n$, there is a vector $X\in T_xM^n$ such that $Z=0$.

By the first variation of $Z$ in $X$
$$\nabla_iR+\frac{1}{n-1}R_{ij}X^j=0, \eqno (2.4)$$
which implies that such a null vector $X$ is unique at each point
and varies smoothly in space-time.

Combining (2.3) and (2.4) we obtain that
$$\frac{\partial R}{\partial t}+\frac{R}{t}+\frac{1}{2}\nabla_iR\cdot X^i=0. \eqno (2.5)$$
By (2.4) and (2.5) and a direct computation, we have
$$\arraycolsep=1.5pt\begin{array}{rcl}
&&X^i(\frac{\partial}{\partial t}-(n-1)\triangle)(\nabla_iR)+
\frac{1}{2(n-1)}X^iX^j(\frac{\partial}{\partial
t}-(n-1)\triangle)R_{ij}\\[4mm]
&&\hskip
0.1cm-\nabla_kR_{ij}\nabla_kX^jX^i-(n-1)\nabla_k\nabla_iR\cdot\nabla_kX^i\\[4mm]
&&\hskip 0.1cm+(\frac{\partial}{\partial
t}-(n-1)\triangle)(\frac{\partial R} {\partial t}+\frac{R}{t})=0,
\end{array}\eqno (2.6)$$
$$\arraycolsep=1.5pt\begin{array}{rcl}
&&(\frac{\partial}{\partial
t}-(n-1)\triangle)(\nabla_iR)=\nabla_i[(\frac{\partial}{\partial
t}-(n-1)\triangle)R]-(n-1)R_{il}\nabla_lR\\[4mm]
&&\hskip 4.1cm=\nabla_i(R^2)-(n-1)R_{il}\nabla_lR,
\end{array}\eqno (2.7)$$
$$\arraycolsep=1.5pt\begin{array}{rcl}
&&(\frac{\partial}{\partial t}-(n-1)\triangle)(\frac{\partial R}
{\partial t}+\frac{R}{t})=3(n-1)R\triangle
R+\frac{1}{2}(n-1)(2-n)|\nabla
R|^2\\[4mm]
&&\hskip 5cm+2R^3+\frac{R^2}{t}-\frac{R}{t^2},
\end{array}\eqno (2.8)$$
$$(\frac{\partial}{\partial
t}-(n-1)\triangle)R_{ij}=\frac{1}{n-2}B_{ij},\eqno (2.9) $$ where
$B_{ij}=(n-1)|Ric|^2g_{ij}+nRR_{ij}-n(n-1)R_{ij}^2-R^2g_{ij}.$ The
combination of (2.6)-(2.9) gives
$$\arraycolsep=1.5pt\begin{array}{rcl}
&&-R(R+\frac{1}{t})^2+\frac{1}{2(n-1)(n-2)}B_{ij}X^iX^j-\frac{1}{2(n-1)}RR_{ij}X^iX^j\\[4mm]
&&\hskip
0.3cm+\frac{n}{2(n-1)}R_{il}R_{jl}X^iX^j+R_{ij}\nabla_kX^i\nabla_kX^j=0.
\end{array}\eqno (2.10)$$

On the other hand, by (2.4) we have
$$\nabla_k\nabla_iR=-\frac{1}{n-1}(X^j\cdot
 \nabla_kR_{ij}+R_{ij}\cdot \nabla_kX^j), \eqno (2.11)$$
and then by taking trace and using the evolution equation of the
scalar curvature,
$$R_{ij}((R+\frac{1}{t})g_{ij}-\nabla_iX^j)=0. \eqno (2.12)$$
Hence it follows from (2.10) and (2.12) that:
$$R_{ij}(\nabla_kX^i-(R+\frac{1}{t})g_{ik})
(\nabla_kX^j-(R+\frac{1}{t})g_{jk})+A_{ij}X^iX^j=0,\eqno (2.13)$$
where
$A_{ij}=\frac{1}{2(n-1)(n-2)}B_{ij}+\frac{1}{2(n-1)}(nR_{il}R_{jl}-RR_{ij}).$

In local coordinate $\{x^i\}$ where $g_{ij}=\delta _{ij}$ and the
Ricci tensor is diagonal, i.e.,
$Ric=diag(\lambda_1,\lambda_2,\cdots,\lambda_n)$, with
$\lambda_1\leq \lambda_2\leq\cdots\leq\lambda_n$, and $e_i,(1\leq
i\leq n)$ is the direction corresponding to the eigenvalue
$\lambda_i$ of the Ricci tensor, we have
$$\sum_i\lambda_i(\nabla_kX^i-(R+\frac{1}{t})g_{ik})^2+A_{ij}X^iX^j=0$$
and $$A_{ij}=diag(\nu_1,\nu_2,\cdots, \nu_n),$$ where
 $$\nu_i=\frac{1}{2(n-1)(n-2)}\sum_{k,l\neq i,k>l}(\lambda_k-\lambda_l)^2\geq 0.$$
So $$\nabla_jX^i=(R+\frac{1}{t})g_{ij},\qquad \mbox{and} \qquad
A_{ij}X^iX^j=0.$$ Thus $\nabla_jX^i$ is symmetric and by the
simply connectedness of $M^n$, there exists a function $f$ such
that $$\nabla_iX^j=\nabla_i\nabla_jf.$$ Hence
$$(R+\frac{1}{t})g_{ij}=\nabla_i\nabla_jf.$$
This means that $g_{ij}(t)$ is a homothetically expanding gradient
soliton.

So we have proved that if the solution exists on $0<t<+\infty$,
and the Harnack quantity $$Z=\frac{\partial R}{\partial t}+\langle
\nabla R,X\rangle+\frac{1}{2(n-1)}R_{ij}X^iX^j+\frac{R}{t}$$
vanishes, then it must be an expanding gradient soliton. If we
have a solution on $\alpha <t<+\infty$, we can replace $t$ by
$t-\alpha$ in the Harnack quantity. Then if $\alpha\rightarrow
-\infty$, the expression $\frac{1}{t-\alpha}\rightarrow 0$ and
disappears. So the Harnack quantity becomes $$Z=\frac{\partial
R}{\partial t}+\langle \nabla
R,X\rangle+\frac{1}{2(n-1)}R_{ij}X^iX^j.$$ Then the rest of the
arguments for the proof of (i) follows.

Hence we complete our proof of Theorem 2.5. $$\eqno \#$$

In order to prove our Theorem 1.1, we need to get more information
about the limit solutions of the Yamabe flow under our
assumptions. So we give two Propositions which are necessary in
our proof in the following section. We first deal with the case of
the Type III limit solutions.

\noindent {\bf Proposition 2.6} \emph{ There exists no noncompact
locally conformally flat Type III limit solution of the Yamabe
flow which satisfies the Ricci pinching condition: $$R_{ij}\geq
\varepsilon\cdot scal \cdot g_{ij}>0,$$ for some $\varepsilon >0$.
}
 \vskip 0.1cm \noindent{\bf Proof.} We argue by contradiction.
 Suppose there is a noncompact locally
 conformally flat Type III limit solution $g_{ij}(t)$ on $M$ which
 satisfies the above Ricci pinching condition. By
 Theorem 2.5, we know that the solution must be a homothetically
 expanding gradient soliton. This means that for any fixed time
 $t=t_0$, we have :
 $$(R+\rho)g_{ij}=\nabla_i\nabla_jf \eqno (2.14)$$
for some positive constant $\rho$ and some function $f$ on $M$.

Differentiating the equation (2.14) and switching the order of
differentiations and then taking trace, we have
$$-(n-1)\nabla_iR=R_{ij}\nabla_jf. \eqno (2.15)$$

Fix the time $t=t_0$ and consider a long shortest geodesic
$\gamma(s)$, $0\leq s\leq \bar{s}$. Let $x_0=\gamma(0)$ and
$X(s)=\dot{\gamma}(s)$. Following by the same arguments as in the
proof of Lemma 1.2 of Perelman \cite{P2} (or see the proof of
Lemma 6.4.1 of \cite{CaZ} for the details) and using the Ricci
pinching condition, we can obtain that
$$|\frac{df}{ds}-\rho s|\leq const.  \eqno(2.16)$$ and
 $$|f-\frac{1}{2}\rho s^2|\leq const\cdot(s+1) \eqno(2.17)$$
for $s$ large enough. From (2.16) and (2.17) we obtain that
$$|\nabla f|^2(x)\geq c\rho f(x)\geq \frac{c}{2}\rho^2s^2=\frac{c}{2}\rho^2d^2(x,x_0)$$
for some constant $c>0$. Then by the same argument as in Theorem I
in \cite{CZ}, we can obtain a contradiction!

Hence we complete the proof of Proposition 2.6. $$\eqno \#$$

For the case of Type II limit solution of the Yamabe flow, we have
the following result:

\noindent {\bf Proposition 2.7} \emph{ Suppose $(M^n,g_{ij}(t))$
is an $n$-dimensional $(n\geq 3)$ complete noncompact locally
conformally flat steady gradient soliton with bounded and positive
Ricci curvature. Assume the scalar curvature assumes its maximum
at a point $p\in M$, then the asymptotic scalar curvature ratio is
infinite, i.e., $$A=\limsup_{s\rightarrow +\infty}Rs^2=+\infty$$
where $s$ is the distance to the point $p$.} \vskip 0.1cm
\noindent{\bf Proof.} We argue by contradiction. Suppose $R\leq
\frac{C}{s^2}$, for some constant $C>0$. By the equation of steady
gradient soliton, we have
$$Rg_{ij}=\nabla_i\nabla_jf, \eqno (2.18)$$ for some smooth
function $f$ on $M$.

Consider the integral curve $\gamma(s), 0\leq s\leq \bar{s}$, of
$\nabla f$ with $\gamma(0)=p$ and $X(s)=\dot{\gamma}(s)$. We first
claim that $M$ is diffeomorphic to $R^n$. Indeed, by
differentiating the equation (2.18) and switching the order of
differentiations and then taking trace, we have
$$-(n-1)\nabla_iR=R_{ij}\nabla_jf. \eqno (2.19)$$
By the positivity of the Ricci curvature, we have
$$(n-1)\nabla_XR+CR\nabla_Xf\geq 0, $$ for some positive constant
$C$ depends only on $n$. This is equivalent to
$$\nabla_X((n-1)\log R+Cf)\geq 0. $$ That is the function $(n-1)\log
R+Cf$ is nondecreasing along $\gamma(s)$.

But by the assumption $$R\leq \frac{C}{s^2},$$ we have $$\log
R\rightarrow -\infty \quad as \quad s\rightarrow +\infty.$$ So
$f(\gamma(s))\rightarrow +\infty $ as $s\rightarrow +\infty.$ That
is $f$ is a exhaustion function on $M$. By (2.18) we know that $f$
is a strictly convex function, so any two level sets of $f$ are
diffeomorphic via the gradient curves of $f$. Combining these and
$f$ is a exhaustion function, we know that $M$  is diffeomorphic
to $R^n$. So we have proved the claim. (We can have another proof
by using the main result of \cite{CH}.)

Next, we follow the argument of Hamilton \cite{Ha95F} to prove
that we can take a limit on $M-\{p\}$ of $g_{ij}(x,t)$ as
$t\rightarrow-\infty$ and the limit is flat.

By (2.18) we have $$\nabla_X\nabla_Xf=R.$$ Integrating it we
obtain
$$X(f(\gamma(s)))-X(f(\gamma(0)))=\int_0^sRds\geq C_0>0$$ for some
constant $C_0>0.$ So we have $|\nabla f|\geq C_0>0$. Then we can
evolve the function $f$ backward with time along the gradient of
$f$. When we go backward in time, this is equivalent to following
outwards along the gradient of $f$, and the speed $|\nabla f|\geq
C_0>0$. So we have
$$\frac{s}{|t|}\geq C_0 \quad as \quad |t| \quad large.$$ Then
$$R\leq \frac{C}{s^2}\leq \frac{C}{C_0^2|t|^2}\quad as \quad |t| \quad large.\eqno(2.20)$$
By the equation of the Yamabe flow, we obtain $$0\geq
\frac{\partial}{\partial t}g_{ij}=-Rg_{ij}\geq
-\frac{C}{C_0^2|t|^2}g_{ij}.$$ Then by the same argument as in
\cite{Ha95F}, we can take a limit on $M-\{p\}$ of $g_{ij}(x,t)$ as
$t\rightarrow-\infty$ and the limit is flat.

Since $M$ is diffeomorphic to $R^n$, we know that $M-\{p\}$ is
diffeomorphic to $S^{n-1}\times R$, but for $n\geq 3$, there
exists no flat metric on it. So we obtain a contradiction.

Hence we complete the proof of the Proposition 2.7.$$\eqno \#$$

\begin{center}{ \bf \large \bf {3. The Proof of the Main Theorem}}\end{center}
 \vskip 0.1cm \noindent{\bf Proof of the Main Theorem 1.1.} We will
 argue by contradiction to prove our Theorem. Let $M^n$ be
 a noncompact conformally flat manifold with nonnegative sectional curvature.
  Suppose $M^n$ has positive and bounded
 scalar curvature and satisfies the Ricci pinching condition:
 $$R_{ij}\geq \varepsilon\cdot scal\cdot g_{ij}$$
 for some $\varepsilon >0$. We evolve the metric by the Yamabe
 flow:
$$
      \left\{
       \begin{array}{lll}
\frac{\partial g_{ij}(x,t)}{\partial t}=-Rg_{ij}(x,t),
          \\[4mm]
  g_{ij}(x,0)=g_{ij}(x),
       \end{array}
    \right.
  \eqno(3.1)
$$
Then by Theorem 2.3 in \cite{CZ02}, we know that the equation has
a smooth solution on a maximal time interval $[0,T)$ with $T>0$
such that either $T=+\infty$ or the evolving metric contracts to a
point at a finite time $T$.

Moreover, for locally conformally flat manifolds, we have
$$R_{ijkl}=\frac{1}{n-2}(R_{ik}g_{jl}+R_{jl}g_{ik}-R_{il}g_{jk}-R_{jk}g_{il})
-\frac{R}{(n-1)(n-2)}(g_{ik}g_{jl}-g_{il}g_{jk}).$$ Then by direct
computation, we have the following evolution equation:
$$\arraycolsep=1.5pt\begin{array}{rcl}
&&\frac{\partial }{\partial t}R_{ijkl}=(n-1)\triangle
R_{ijkl}-R\cdot R_{ijkl}
+\frac{n-1}{n-2}[(R_{imkn}R_{mn}-R_{ik}^2)g_{jl}\\[4mm]
&&\hskip 1.7cm+(R_{jmln}R_{mn}-R_{jl}^2)g_{ik}-(R_{jmkn}R_{mn}-R_{jk}^2)g_{il}\\[4mm]
&&\hskip 1.7cm-(R_{imln}R_{mn}-R_{il}^2)g_{jk}]\\[4mm]
&&\hskip 1.4cm=(n-1)\triangle R_{ijkl}-R\cdot
R_{ijkl}+\frac{1}{(n-2)^2}(B_{ik}g_{jl}+B_{jl}g_{ik}-B_{il}g_{jk}-B_{jk}g_{il}),
\end{array}$$
where
$B_{ij}=(n-1)|Ric|^2g_{ij}+nRR_{ij}-n(n-1)R_{ij}^2-R^2g_{ij}.$ In
a moving frame, we have:
$$\arraycolsep=1.5pt\begin{array}{rcl}
&&\frac{\partial }{\partial t}R_{abcd}=(n-1)\triangle
R_{abcd}-R\cdot
R_{abcd}+\frac{1}{(n-2)^2}(B_{ac}g_{bd}+B_{bd}g_{ac}-B_{ad}g_{bc}-B_{bc}g_{ad})\\[4mm]
&&\hskip
1.7cm+\frac{R}{n-2}\cdot(R_{ac}g_{bd}+R_{bd}g_{ac}-R_{ad}g_{bc}-R_{bc}g_{ad}).
\end{array}$$
At a point where $g_{ab}=\delta_{ab}$ and the Ricci tensor is
diagonal:
$$Ric=diag(\lambda_1,\lambda_2,\cdots,\lambda_n),$$ with
$\lambda_1\leq \lambda_2\leq\cdots\leq\lambda_n$, we also have
$B_{ab}$ is diagonal and the sectional curvature
$$R_{abab}=\frac{1}{n-2}(\lambda_a+\lambda_b)-\frac{R}{(n-1)(n-2)}.$$
If at some point, the sectional curvature $R_{1212}=0$, then
$\lambda_1+\lambda_2=\frac{R}{n-1}$. Hence if $n\geq 4$, we have :
$$\arraycolsep=1.5pt\begin{array}{rcl}
&&\hskip
0.5cm\frac{1}{(n-2)^2}(B_{aa}+B_{bb})+\frac{R}{n-2}(\lambda_a+\lambda_b)\\[4mm]
&&\hskip
0.1cm=\frac{1}{(n-2)^2}[2(n-1)|Ric|^2+nR(\lambda_a+\lambda_b)
-n(n-1)(\lambda_a^2+\lambda_b^2)-2R^2]+\frac{R^2}{(n-2)(n-1)}\\[4mm]
&&\hskip 0.1cm\geq
\frac{1}{(n-2)^2}[\frac{2(n-1)}{n}R^2+\frac{nR^2}{n-1}-n(n-1)\frac{R^2}
{(n-1)^2}-2R^2]+\frac{R^2}{(n-2)(n-1)}\\[4mm]
&&\hskip 0.1cm=\frac{n^2-4n+2}{n(n-1)(n-2)^2}R^2\\[4mm]
&&\hskip 0.1cm>0,
\end{array}$$
if $n=3$, by direct calculation, we have:
$$\arraycolsep=1.5pt\begin{array}{rcl}
&&\hskip
0.5cm\frac{1}{(n-2)^2}(B_{11}+B_{22})+\frac{R}{n-2}(\lambda_1+\lambda_2)\\[4mm]
&&\hskip 0.1cm=B_{11}+B_{22}+\frac{1}{2}R^2\\[4mm]
&&\hskip
0.1cm=4|Ric|^2+3R(\lambda_1+\lambda_2)-6(\lambda_1^2+\lambda_2^2)-2R^2+\frac{1}{2}R^2\\[4mm]
&&\hskip
0.1cm=4(\lambda_1^2+\lambda_2^2+\lambda_3^2)-6(\lambda_1^2+\lambda_2^2)\\[4mm]
&&\hskip 0.1cm=4\lambda_3^2-2(\lambda_1^2+\lambda_2^2)\\[4mm]
&&\hskip 0.1cm=R^2-2(\lambda_1^2+\lambda_2^2)\\[4mm]
&&\hskip
0.1cm=\lambda_1^2+\lambda_2^2+\lambda_3^2+2\lambda_1\lambda_2+2(\lambda_1+\lambda_2)\lambda_3
-2(\lambda_1^2+\lambda_2^2)\\[4mm]
&&\hskip
0.1cm=(\lambda_3^2-\lambda_1^2)+(\lambda_3^2-\lambda_2^2)+2\lambda_1\lambda_2\\[4mm]
&&\hskip 0.1cm>0.
\end{array}$$
So we obtain that the nonnegative sectional curvature is preserved
under the Yamabe flow.

Next we claim that under our assumption, the solution $g_{ij}(t)$
has a long-time existence. Otherwise, using the same argument as
in Theorem 1.2 in \cite{Chow2}, we know that the Ricci pinching
condition is preserved under the Yamabe flow. Then by a scaling
argument as in Ricci flow, we can take a limit to obtain a
noncompact solution to the Yamabe flow with constant positive
Ricci curvature, which is a contradiction with Bonnet-Myers'
Theorem. So we have the long-time existence result.

By a standard rescaling argument similarly as in Ricci flow, we
know that there exists a sequence of dilations of the solution
which converges to a noncompact limit solution, which we also
denote by $g_{ij}(t)$, of Type II or Type III with positive scalar
curvature and it still satisfies the Ricci pinching condition.

Now we consider its universal covering space, then we also have a
solution on its universal cover which is of Type II or Type III.
So in the following we consider the limit solution is defined on
its universal cover.

If the limit solution is of Type III, then by Theorem 2.5, we know
that it is a homothetically expanding gradient soliton, but from
Proposition 2.6, we know that there exists no such limit solution
of Type III satisfies the Ricci pinching condition. So the limit
must be of Type II.

Suppose the limit solution is of Type II, then by Theorem 2.5, we
know that it is a homothetically steady gradient soliton. From
Proposition 2.7, we also know that $$\limsup_{s\rightarrow
+\infty}Rs^2=+\infty,$$ where $s$ is the distance function from
the point $p$ where the scalar curvature $R$ assumes its maximum.
Then by the result of Hamilton \cite{Ha95F}, we can take a
sequence of points $x_k$ divergent to infinity and a sequence of
$r_k$, such that $r_k^2R(x_k)\rightarrow+\infty$ and
$\frac{d(p,x_k)}{r_k}\rightarrow+\infty$ and $$R(x)\leq 2 R(x_k)$$
for all points $x\in B(x_k,r_k)$. Then by a same argument as in
Ricci flow, we obtain that $(M,R(x_k)g_{ij},x_k)$ converge to a
limit manifold $(\widetilde{M},\widetilde{g_{ij}},\widetilde{x})$
with nonnegative sectional curvature. By Proposition 2.3 in
\cite{CZ06}, we know that the limit manifold will split a line.
Since the Ricci pinching condition is preserved under dilations,
we conclude that the limit must be also satisfies the Ricci
pinching condition. And this is a contradiction.

Therefore the proof of the main theorem 1.1 is completed.$$\eqno
\#$$


\begin{thebibliography}{99}

\bibitem{BW}   C. B$\ddot{o}$hm, and B. Wilking,  {\it Manifolds with positive curvature
 operators are space forms}, arXiv:math.DG/0606187 June 2006.

\bibitem{Cao}  Cao, H. D.,  {\it Limits of solutions to the
K$\ddot{a}$hler Ricci flow},   J. Diff. Geom. {\bf 45}, (1997),
257-272.

\bibitem{CaZ}  Cao, H. D. and  Zhu, X. P.,  {\it A complete proof of the
Poincar$\acute{e}$ and geometrization conjecture -- application of
the Hamilton-Perelman theory of the Ricci flow,} Asian J. Math. {\bf
10} (2006), no. 2, 165-492.

\bibitem{CH}  G. Carron, and M. Herzlich, {\it Conformally  flat
manifolds with nonnegative Ricci curvature}, arXiv:math.DG/0410180
v2 Oct 2004.

\bibitem{CZ} Chen, B. L. and Zhu, X. P., {\it Complete
Riemannian manifolds with pointwise pinched curvature}, Invent.
Math. {\bf 140} (2000) no. 2, 423-452.

\bibitem{CZ02} Chen, B. L. and Zhu, X. P., {\it A gap theorem for
 complete noncompact manifolds with nonnegative Ricci curvature },
 Comm. Anal. Geom.  {\bf 10} (2002)  no.1,   217-239.

\bibitem{CZ06} Chen, B. L. and Zhu, X. P., {\it  Ricci flow with surgery on four-manifolds
with positive isotropic curvature},  J. Diff. Geom. {\bf 74},
(2006) no.2, 177-264.



\bibitem{Chow2} Chow, B., {\it The Yamabe flow on locally conformally flat
manifolds with positive Ricci curvature },  Comm. Pure Appl.
Math., {\bf XIV}  (1992),   1003-1014.


\bibitem{Ha93} Hamilton, R. S., {\it The Harnack estimate for the Ricci
flow},  J. Diff. Geom. {\bf 37} (1993), 225-243.

\bibitem{Ha93E} Hamilton, R. S., {\it Eternal solutions to the
Ricci flow},   J. Diff. Geom. {\bf 38} (1993), 1-11.

\bibitem{Ha94} Hamilton, R. S., {\it  Convex hypersurfaces with
pinched second fundamental form},  Comm. Anal. Geom.  {\bf 2}
(1994),  no.1,  167-172.

\bibitem{Ha95F}  Hamilton, R. S., {\it The formation of singularities in
the Ricci flow}, Surveys in Differential Geometry (Cambridge, MA,
1993), {\bf 2}, 7-136, International Press, Combridge, MA,1995.


\bibitem{NW}  Ni, L. and Wu, B. Q., {\it Complete manifolds with
nonnegative curvature operator}, arXiv:math.DG/0607356  v1 July
2006.

\bibitem{P2}G. Perelman, {\it Ricci flow with surgery on three
manifolds},  arXiv:math. DG/0303109  v1 March 10, 2003. prepeint.


\end{thebibliography}
\end{document}